\documentclass[12pt,reqno]{amsart}
\usepackage{amssymb,verbatim,enumerate,ifthen}
\usepackage[mathscr]{eucal}
\usepackage[utf8]{inputenc}
\usepackage[T1]{fontenc}
\allowdisplaybreaks[1]
\def\N{\mathbb{N}}
\def\R{\mathbb{R}}

\newtheorem{theorem}{Theorem}
\newtheorem*{theorem*}{Theorem}
\def\Thm#1#2{\ifthenelse{\equal{#1}{*}}{\begin{theorem*}#2\end{theorem*}}
             {\begin{theorem}\label{T#1}#2\end{theorem}}}
\newtheorem{Atheorem}{Theorem}

\def\thm#1{Theorem~\ref{T#1}}

\newtheorem{proposition}[theorem]{Proposition}
\newtheorem*{proposition*}{Proposition}
\def\Prp#1#2{\ifthenelse{\equal{#1}{*}}{\begin{proposition*}#2\end{proposition*}}
             {\begin{proposition}\label{P#1}#2\end{proposition}}}
\def\prp#1{Proposition~\ref{P#1}}

\newtheorem{corollary}[theorem]{Corollary}
\newtheorem*{corollary*}{Corollary}
\def\Cor#1#2{\ifthenelse{\equal{#1}{*}}{\begin{corollary*}#2\end{corollary*}}
             {\begin{corollary}\label{C#1}#2\end{corollary}}}
\def\cor#1{Corollary~\ref{C#1}}

\newtheorem{lemma}[theorem]{Lemma}
\newtheorem*{lemma*}{Lemma}
\def\Lem#1#2{\ifthenelse{\equal{#1}{*}}{\begin{lemma*}#2\end{lemma*}}
             {\begin{lemma}\label{L#1}#2\end{lemma}}}
\def\lem#1{Lemma~\ref{L#1}}
\newtheorem{Alemma}{Lemma}

\theoremstyle{definition}
\newtheorem{remark}[theorem]{Remark}
\newtheorem*{remark*}{Remark}
\def\Rem#1#2{\ifthenelse{\equal{#1}{*}}{\begin{remark*}\rm #2\end{remark*}}
             {\begin{remark}\label{R#1}\rm #2\end{remark}}}

\newtheorem{example}[theorem]{Example}
\newtheorem*{example*}{Example}
\def\Exa#1#2{\ifthenelse{\equal{#1}{*}}{\begin{example*}\rm #2\end{example*}}
             {\begin{example}\label{Ex#1}\rm #2\end{example}}}

\def\eq#1{{\rm(\ref{E#1})}}
\def\Eq#1#2{\ifthenelse{\equal{#1}{*}}
  {\begin{equation*}\begin{aligned}#2\end{aligned}\end{equation*}}
  {\begin{equation}\begin{aligned}\label{E#1}#2\end{aligned}\end{equation}}}

\def\diag{\mathop{\hbox{\rm diag}}\nolimits}
\def\conv{\mathop{\hbox{\rm conv}}\nolimits}
\def\id{\mathop{\text{\rm id}}\nolimits}

\begin{document}
\vspace{5mm}

\date{\today}

\title[Equality of generalized Bajraktarević means]{On the equality of generalized Bajraktarević means under first-order differentiability assumptions}

\author[Zs. P\'ales]{Zsolt P\'ales}
\address{Institute of Mathematics, University of Debrecen, 
H-4002 Debrecen, Pf. 400, Hungary}
\email{pales@science.unideb.hu}

\author[A. Zakaria]{Amr Zakaria}
\address[A. Zakaria]{Department of Mathematics, Faculty of Education, Ain Shams University, Cairo 11341, Egypt}
\email{amr.zakaria@edu.asu.edu.eg}

\thanks{The research of the second author was supported by the Hungarian Scientific Research Fund (OTKA) Grant K- 134191.}

\dedicatory{Dedicated to the 100th birthday of Professor János Aczél}

\begin{abstract}
In this paper we consider the equality problem of generalized Bajraktarević means, i.e., we are going to solve the functional equation
\begin{equation}\label{E0}\tag{*}
	f^{(-1)}\bigg(\frac{p_1(x_1)f(x_1)+\dots+p_n(x_n)f(x_n)}{p_1(x_1)+\dots+p_n(x_n)}\bigg)=g^{(-1)}\bigg(\frac{q_1(x_1)g(x_1)+\dots+q_n(x_n)g(x_n)}{q_1(x_1)+\dots+q_n(x_n)}\bigg),
\end{equation}
which holds for all $x=(x_1,\dots,x_n)\in I^n$, where $n\geq 2$, $I$ is a nonempty open real interval, the unknown functions $f,g:I\to\R$ are strictly monotone, $f^{(-1)}$ and $g^{(-1)}$ denote their generalized left inverses, respectively, and the vector-valued weight functions $p=(p_1,\dots,p_n):I\to\R_{+}^n$ and $q=(q_1,\dots,q_n):I\to\R_{+}^n$ are also unknown. This equality problem in the symmetric two-variable case (i.e., when $n=2$ and $p_1=p_2$, $q_1=q_2$) was solved under sixth-order regularity assumptions by Losonczi in 1999. The authors of this paper improved this result in 2023 by reaching the same conclusion assuming only first-order differentiability. In the nonsymmetric case, assuming third-order differentiability of $f$, $g$ and the first-order differentiability of at least three of the functions $p_1,\dots,p_n$, Grünwald and Páles proved that \eq{0} holds if and only if there exist four constants $a,b,c,d\in\R$ with $ad\neq bc$ such that 
\Eq{*}{
        cf+d>0,\qquad 
	g=\frac{af+b}{cf+d},\qquad\mbox{and}\qquad q_\ell=(cf+d)p_\ell\qquad (\ell\in\{1,\dots,n\}).
}
The main goal of this paper is to establish the same conclusion under first-order differentiability. 
\end{abstract}

\maketitle

\section{Introduction}

Throughout this paper, the symbols $\R$ and $\R_+$ will stand for the sets of real and positive real numbers, respectively, and $I$ will always denote a nonempty open real interval. In theory of quasi-arithmetic means the characterization of the equality of means with different generators is a basic problem which was completely solved in the book \cite{HarLitPol34}. Using this characterization, the homogeneous quasi-arithmetic means can also be found: they are exactly the power means and the geometric mean. In \cite{Baj58} (cf.\ also \cite{Baj63}) Bajraktarević introduced a new generalization of quasi-arithmetic means by adding a weight function to the formula of quasi-arithmetic means. He also characterized the equality of such means (called Bajraktarević means since then) in the at least $3$-variable setting assuming three times differentiability. Daróczy and Losonczi in \cite{DarLos70}, later Daróczy and Páles in \cite{DarPal82} arrived at the same conclusion with first-order differentiability and without differentiability, respectively, but assuming the equality for all $n\in\N$. As an application of the characterization of the equality, Aczél and Daróczy in \cite{AczDar63c} determined the homogeneous Bajraktarević means that include the means which were introduced by Gini in \cite{Gin38}. Losonczi in \cite{Los99} described the equality of two-variable Bajraktarević means under sixth-order regularity assumptions and an algebraic condition which was later removed in \cite{Los06b}. A new form of the result of Losonczi was established in the paper \cite{LosPalZak21}. Using these results, the homogeneous two-variable means were also determined by Losonczi \cite{Los07a}, \cite{Los07b}. 

A characterization of the equality of two-variable Bajraktarević means has been proved in \cite[Theorem 7 and Corollary 8]{PalZak23} by the authors assuming only first-order differentiability of the generating functions of the means.

The purpose of this paper is to consider an analogous problem for an nonsymmetric extension of Bajraktarević means, which was introduced by \cite{GruPal20}. In Theorem 12 of this paper, the equality problem of the extended means was solved under three times differentiability conditions. The main result of this paper will establish the same conclusion under first-order differentiability only.

Equality problems related to other generalizations of Bajraktarević means were investigated and solved in the papers \cite{PalZak20b,PalZak20a,PalZak22a}.

\section{Preliminary and auxiliary results}

Given a subset $S\subseteq\R$, the smallest convex set containing $S$, which is identical to the smallest interval containing $S$, will be denoted by $\conv(S)$.
For our definition of generalized Bajraktarević means, we shall need the following lemma about the existence and properties of the left inverse of strictly monotone (but not necessarily continuous) functions.

\Lem{SMF}{{\rm(\cite[Lemma 1]{GruPal20})}
Let $f:I\to\R$ be a strictly monotone function. Then there exists a uniquely determined monotone function $f^{(-1)}:\conv(f(I))\to I$ such that $f^{(-1)}$ is the left inverse of $f$, i.e.,
\Eq{SMF1}{
   (f^{(-1)}\circ f)(x)=x \qquad(x\in I).
}
Furthermore, $f^{(-1)}$ is monotone in the same sense as $f$, continuous,
\Eq{SMF2}{
  (f\circ f^{(-1)})(y)=y \qquad(y\in f(I)),
}
and
\Eq{SMF3}{
  \liminf_{x\to f^{(-1)}(y)}f(x)\leq y \leq \limsup_{x\to f^{(-1)}(y)}f(x) 
  \qquad(y\in\conv(f(I))).
}
Thus, if $f$ is lower (resp.\ upper) semicontinuous at $f^{(-1)}(y)$, then $f\circ f^{(-1)}(y)\leq y$ (resp.\ $y\leq f\circ f^{(-1)}(y)$).}

It is clear from \eq{SMF1} and \eq{SMF2} that the restriction of $f^{(-1)}$ to $f(I)$ is the inverse of $f$ in the standard sense. Therefore, $f^{(-1)}$ is the continuous and monotone extension of the inverse of $f$ to the smallest interval containing the range of $f$.

Given a strictly monotone function $f:I\to\R$ and an $n$-tuple of positive valued functions $p=(p_1,\dots,p_n):I\to\R_{+}^n$, we introduce the \emph{$n$-variable generalized Bajraktarević mean} $A_{f,p}:I^n\to I$ by the following formula:
\Eq{BM}{
	A_{f,p}(x):=f^{(-1)}\bigg(\frac{p_1(x_1)f(x_1)+\dots+p_n(x_n)f(x_n)}{p_1(x_1)+\dots+p_n(x_n)}\bigg) \qquad (x=(x_1,\dots,x_n)\in I^n),
}
and, to simplify the notations, we will use the following definition:
\Eq{Rfp}{
	R_{f,p}(x)
	:=\frac{p_1(x_1)f(x_1)+\dots+p_n(x_n)f(x_n)}{p_1(x_1)+\dots+p_n(x_n)}.
}

\Thm{BM1}{{\rm(\cite[Theorem 2]{GruPal20})}
Let $f:I\to\R$ be strictly monotone and $p=(p_1,\dots,p_n):I\to\R_{+}^n$. Then the function $A_{f,p}:I^n\to I$ given by \eq{BM} is well-defined and it is a mean, that is, 
\Eq{MV}{
	\min(x)\leq A_{f,p}(x)\leq \max(x)\qquad (x=(x_1,\dots,x_n)\in I^n).
}
}

\Thm{BM2}{{\rm(\cite[Theorem 3]{GruPal20})}
Let $f:I\to\R$ be strictly increasing and $p=(p_1,\dots,p_n):I\to\R_{+}^n$. Then, for all $x=(x_1,\dots,x_n)\in I^n$, the equality $y=A_{f,p}(x)$ holds if and only if 
	\Eq{MV2}{
	\sum_{i=1}^{n}p_i(x_i)(f(z)-f(x_i))
	\begin{cases}
	 <0 & \mbox{ for } z\in I,\, z<y,\\[2mm]
	 >0 & \mbox{ for } z\in I,\, z>y.
	\end{cases}
	}
If $f$ is strictly decreasing, then the inequalities \eq{MV2} hold with reversed inequality sign.}

\Cor{BM2}{{\rm(\cite[Corollary 4]{GruPal20})}
Let $f:I\to\R$ be continuous, strictly monotone, and $p=(p_1,\dots,p_n):I\to\R_{+}^n$. Then, for all $x=(x_1,\dots,x_n)\in I^n$, the value $y=A_{f,p}(x)$ is the unique solution of the equation
\Eq{MV3}{
	\sum_{i=1}^{n}p_i(x_i)(f(y)-f(x_i))=0.
}}

For the formulation of the subsequent results, we define the \emph{diagonal} $\diag(I^n)$ of $I^n$ by 
\Eq{*}{
\diag(I^n):=\{(x,\dots,x)\in\R^n\mid x\in I\}
}
and, for a function $F:I^n\to\R$, its \emph{diagonalization $F^\Delta:I\to\R$} is given as 
\Eq{*}{
F^\Delta(x):=F(x,\dots,x) \qquad(x\in I).
}

Given two functions $F,G:I^n\to\R$, we say that \emph{$F$ and $G$ are locally equal on $I^n$} (in other words, \emph{$F$ and $G$ are equal near the diagonal of $I^n$}) if there exists an open set $U\subseteq I^n$ containing $\diag(I^n)$ such that $F(x)=G(x)$ is valid for all $x\in U$. If the equality $F(x)=G(x)$ is valid for all $x\in I^n$, then we say that \emph{$F$ and $G$ are globally equal on $I^n$}.

Given $p=(p_1,\dots,p_n):I\to\R_{+}^n$ and $q=(q_1,\dots,q_n):I\to\R_{+}^n$, we will also use the following notations:
\Eq{*}{
	p_0:=p_1+\dots+p_n,\qquad q_0:=q_1+\dots+q_n,\qquad\mbox{and}\qquad r_0:=\frac{q_0}{p_0}.
} 

\Thm{DB}{{\rm(\cite[Assertion (1) of Theorem 8]{GruPal20})}
Let $f: I\to\R$ be a differentiable function on $I$ with a nonvanishing first derivative, let $p=(p_1,\dots,p_n):I\to\R_{+}^n$ and let $i\in \{1,\dots,n\}$. If $p_i$ is continuous on $I$, then the first-order partial derivative $\partial_i A_{f,p}$ exists on $\diag(I^n)$ and
\Eq{*}{
\big(\partial_i A_{f,p}\big)^\Delta=\frac{p_i}{p_0}.
}}

The next result establishes a sufficient condition for the equality of the $n$-variable generalized Bajraktarević means. We will call this situation the \emph{canonical case of the equality}.

\Thm{BM3}{{\rm(\cite[Theorem 5]{GruPal20})}
Let $f,g:I\to\R$ be strictly monotone and $p,q:I\to\R_{+}^n$. If there exist $a,b,c,d\in\R$ with $ad\neq bc$ such that 
\Eq{fg}{
    cf+d>0,\qquad
	g=\frac{af+b}{cf+d},\qquad\mbox{and}\qquad q=(cf+d)p
}
hold on $I$, then the $n$-variable generalized Bajraktarević means $A_{f,p}$ and $A_{g,q}$ are identical on $I^n$, i.e., globally equal on $I^n$.}

With the aid of the following lemma, we can reduce the regularity assumptions in our statements. 

\Lem{reg}{{\rm(\cite[Lemma 6]{GruPal20})}
Let $f,g:I\to\R$ be continuous strictly monotone functions, $n\geq2$, and $p=(p_1,\dots,p_n):I\to\R_{+}^n,\,q=(q_1,\dots,q_n):I\to\R_{+}^n$. Assume that $A_{f,p}$ and $A_{g,q}$ are locally equal on $I^n$.  Then the following two assertions hold.
\begin{enumerate}[(i)]
\item For all $i\in\{1,\dots,n\}$, the function $p_i$ is continuous on $I$ if and only if the function $q_i$ is continuous on $I$.
\item Assume that $f,g:I\to\R$ are  differentiable (resp.\ continuously differentiable) functions on $I$ with nonvanishing first derivatives. Then, for all $i\in\{1,\dots,n\}$, the function $p_i$ is differentiable (resp.\ continuously differentiable) on $I$ if and only if $q_i$ is differentiable (resp.\ continuously differentiable) on $I$.
\end{enumerate}
}

The following theorem is of basic importance for our investigations. 

\Thm{BM4}{{\rm(\cite[Theorem 7]{GruPal20})}
Let $f,g:I\to\R$ be continuous, strictly monotone, $n\geq2$ and $p:I\to\R_{+}^n$ be a continuous function on $I$. Let further $q:I\to\R_{+}^n$. Assume that $A_{f,p}$ and $A_{g,q}$ are locally equal on $I^n$ and that there exist $a,b,c,d\in\R$ with $ad\neq bc$ and a nonempty open subinterval $I_0$ of $I$ such that \eq{fg} holds on $I_0$. Then \eq{fg} is also valid on $I$.}

In the next result we point out that the local equality of the $n$-variable Bajraktarević means implies the local equality of the $k$-variable descendant Bajraktarević means for $k\leq n$.

\Thm{RV}{Let $f,g:I\to\R$ be strictly monotone and $p=(p_1,\dots,p_n)$, $q=(q_1,\dots,q_n):I\to\R_{+}^n$. Assume that the $n$-variable $A_{f,p}$ and $A_{g,q}$ are locally (globally) equal on $I^n$. Then, for any $k\in\N$ and $i_1,\dots,i_k\in\{1,\dots,n\}$ with $i_1<\dots<i_k$, we have that the $k$-variable means $A_{f,(p_{i_1},\dots,p_{i_k})}$ and $A_{q,(q_{i_1},\dots,q_{i_k})}$ are locally (globally) equal on $I^k$.}

\begin{proof} Without loss of generality, we may assume that $f$ and $g$ are strictly increasing functions. Then $f^{(-1)}$ and $g^{(-1)}$ are continuous increasing functions. 

According to the local equality of $A_{f,p}$ and $A_{g,q}$ on $I^n$, there exists an open set $U\subseteq I^n$ containing the diagonal of $I^n$ such that, for all $x\in U$, we have $A_{f,p}(x)=A_{g,q}(x)$. 

In what follows, for all $i\in\{1,\dots,n\}$, let $e_i\in\R^n$ denote the $i$th vector of the standard base of $\R^n$, i.e., let $e_i:=(\delta_{i,j})_{j=1}^n$, where $\delta_{i,j}$ stands for the Kronecker symbol.

Let us now fix $k\in\N$ and $i_1,\dots,i_k\in\{1,\dots,n\}$ with $i_1<\dots<i_k$ and define the set $V\subseteq I^k$ by
\Eq{*}{
  V:=
  \bigg\{(y_1,\dots,y_k)\in I^k\,\bigg|\,\forall \,t\in\!\Big[\min_{1\leq\alpha\leq k}y_\alpha,\min_{1\leq\alpha\leq k}y_\alpha\Big]\colon \,\sum_{j=1}^n \Big[t+\sum_{\alpha=1}^k\delta_{j,i_\alpha}(y_\alpha-t)\Big]e_j\in U\bigg\}.
}
The interior of $V$ (denoted as $V^\circ$) is an open subset of $I^k$. To see that $V^\circ$ contains the diagonal of $I^k$, let $y\in I$ be fixed. Then the vector 
\Eq{*}{
\sum_{j=1}^n \big[t+\delta_{j,i_1}(y_1-t)+\dots+\delta_{j,i_k}(y_k-t)\big]e_j
}
depends continuously on $(y_1,\dots,y_k,t)$ and it
tends to the diagonal vector $(y,\dots,y)\in U\subseteq I^n$ if the variables $y_1,\dots,y_k$, and $t$ tend to $y$. In other words, due to the openness of $U$, there exists a positive number $r$ such that
\Eq{*}{
  \sum_{j=1}^n \big[t+\delta_{j,i_1}(y_1-t)+\dots+\delta_{j,i_k}(y_k-t)\big]e_j\in U 
  \qquad\mbox{for all}\qquad
  (y_1,\dots,y_k,t)\in (y-r,y+r)^{k+1}.
}
In particular, this inclusion holds if $y-r<\min(y_1,\dots,y_k)\leq t\leq \max(y_1,\dots,y_k)<y+r$, hence, in this case $(y_1,\dots,y_k)\in V$. This shows that $(y,\dots,y)\in I^k$ is an interior point of $V$ and hence $V^\circ$ is an open subset of $I^k$ which contains the diagonal of $I^k$. In addition, if $A_{f,p}$ and $A_{g,q}$ are globally equal on $I^n$, then $U=I^n$ and $V=V^\circ=I^k$.

We are now going to show that the means $A_{f,(p_{i_1},\dots,p_{i_k})}$ and $A_{q,(q_{i_1},\dots,q_{i_k})}$ are equal on $V^\circ$. To the contrary, assume that, for some $(y_1,\dots,y_k)\in V^\circ$, we have
\Eq{*}{
  f^{(-1)}\bigg(\frac{p_{i_1}(y_1)f(y_1)+\dots+p_{i_k}(y_k)f(y_k)}{p_{i_1}(y_1)+\dots+p_{i_k}(y_k)}\bigg) 
  \neq g^{(-1)}\bigg(\frac{q_{i_1}(y_1)g(y_1)+\dots+q_{i_k}(y_k)g(y_k)}{q_{i_1}(y_1)+\dots+q_{i_k}(y_k)}\bigg).
}
We may assume that
\Eq{*}{
  f^{(-1)}\bigg(\frac{p_{i_1}(y_1)f(y_1)+\dots+p_{i_k}(y_k)f(y_k)}{p_{i_1}(y_1)+\dots+p_{i_k}(y_k)}\bigg) 
  <g^{(-1)}\bigg(\frac{q_{i_1}(y_1)g(y_1)+\dots+q_{i_k}(y_k)g(y_k)}{q_{i_1}(y_1)+\dots+q_{i_k}(y_k)}\bigg).
}
The set of points where $f$ or $g$ are discontinuous is countable subset of $I$, therefore, there exists a point $t\in I$, where $f$ and $g$ are continuous and
\Eq{t}{
  f^{(-1)}\bigg(\frac{p_{i_1}(y_1)f(y_1)+\dots+p_{i_k}(y_k)f(y_k)}{p_{i_1}(y_1)+\dots+p_{i_k}(y_k)}\bigg) 
  <t<g^{(-1)}\bigg(\frac{q_{i_1}(y_1)g(y_1)+\dots+q_{i_k}(y_k)g(y_k)}{q_{i_1}(y_1)+\dots+q_{i_k}(y_k)}\bigg).
}
Due to the mean value property of the means (established in \thm{BM1}), it follows that
\Eq{*}{
  \min(y_1,\dots,y_k)<t<\min(y_1,\dots,y_k),
}
therefore, 
\Eq{*}{
  x=(x_1,\dots,x_n):=\sum_{j=1}^n \Big[t+\sum_{\alpha=1}^k\delta_{j,i_\alpha}(y_\alpha-t)\Big]e_j\in U,
}
and hence $A_{f,p}(x)=A_{g,q}(x)$. On the other hand, the inequalities in \eq{t} imply that
\Eq{*}{
  \frac{p_{i_1}(y_1)f(y_1)+\dots+p_{i_k}(y_k)f(y_k)}{p_{i_1}(y_1)+\dots+p_{i_k}(y_k)}
  <f(t),\qquad g(t)<\frac{q_{i_1}(y_1)g(y_1)+\dots+q_{i_k}(y_k)g(y_k)}{q_{i_1}(y_1)+\dots+q_{i_k}(y_k)}.
}
Rearranging these inequalities, we get
\Eq{*}{
  \sum_{\alpha=1}^k p_{i_\alpha}(y_{\alpha})(f(y_\alpha)-f(t))<0,
  \qquad
  \sum_{\alpha=1}^k q_{i_\alpha}(y_{\alpha})(g(y_\alpha)-g(t))>0.
}
Now taking into account the equalities
\Eq{xj}{
   x_j
   :=\begin{cases}
	t & \mbox{if } j\not\in\{i_1,\dots,i_k\},\\[1mm]
	y_{\alpha} & \mbox{if } j=i_\alpha \quad (\alpha\in\{1,\dots,k\}),
   \end{cases}
}
it follows that
\Eq{*}{
  \sum_{i=1}^n p_i(x_i)(f(x_i)-f(t))<0,
  \qquad
  \sum_{i=1}^n q_i(x_i)(g(x_i)-g(t))>0.
}
From these inequalities we obtain that
\Eq{*}{
  \frac{p_{1}(x_1)f(x_1)+\dots+p_{n}(x_n)f(x_n)}{p_{1}(x_1)+\dots+p_{n}(x_n)}
  <f(t),\qquad g(t)<\frac{q_{1}(x_1)g(x_1)+\dots+q_{n}(x_n)g(x_n)}{q_{1}(x_1)+\dots+q_{n}(x_n)}.
}
Finally, using again that $f$ and $g$ are continuous at the point $t$, we can conclude that
\Eq{*}{
  f^{(-1)}\bigg(\frac{p_{1}(x_1)f(x_1)+\dots+p_{n}(x_n)f(x_n)}{p_{1}(x_1)+\dots+p_{n}(x_n)}\bigg)
  <t<g^{(-1)}\bigg(\frac{q_{1}(x_1)g(x_1)+\dots+q_{n}(x_n)g(x_n)}{q_{1}(x_1)+\dots+q_{n}(x_n)}\bigg),
}
which contradicts the equality $A_{f,p}(x)=A_{g,q}(x)$.
\end{proof}

\Cor{RV}{
Let $n\geq2$, $f,g:I\to\R$ be strictly monotone and $p=(p_1,\dots,p_n),q=(q_1,\dots,q_n):I\to\R_{+}^n$.
Assume that the $n$-variable means $A_{f,p}$ and $A_{g,q}$ are locally (globally) equal on $I^n$. Then, for all $i,j\in\{1,\dots,n\}$ with $i<j$, the two-variable means $A_{f,(p_i,p_j)}$ and $A_{g,(q_i,q_j)}$ are locally (globally) equal on $I^2$.}

\section{Reduction of the equality problem}

In our first result, with a certain substitution, we reduce the number of unknown functions in the local (global) equality problem of $n$-variable generalized Bajraktarević means from $2n+2$ to $2n+1$. 

\Thm{Red}{Let $f,g:I\to\R$ be continuous strictly monotone functions, $n\geq 2$ and $p=(p_1,\dots,p_n) :I\to\R_+^n$, $q=(q_1,\dots,q_n):I\to\R_+^n$. Then the means $A_{f,p}$ and $A_{g,q}$ are locally (globally) equal on $I^n$ if and only if, the means 
$A_{\id,P}$ and $A_{h,Q}$ are locally (globally) equal on $J^n$, where the interval $J$, the functions $h:J\to\R$, and $P,Q:J\to\R_+^n$ are defined by  
\Eq{hPQ}{
  J:=f(I),\qquad
  h:=g\circ f^{-1}, \qquad
  P:=p\circ f^{-1},\qquad\mbox{and}\qquad
  Q:=q\circ f^{-1}.
}}

\begin{proof} Assume that the means $A_{f,p}$ and $A_{g,q}$ are locally equal on $I^n$. Then the equality 
\Eq{Eq1}{
  f^{-1}\bigg(\frac{p_1(x_1)f(x_1)+\dots +p_n(x_n)f(x_n)}{p_1(x_1)+\dots+p_n(x_n)}\bigg)
  =g^{-1}\bigg(\frac{q_1(x_1)g(x_1)+\dots +q_n(x_n)g(x_n)}{q_1(x_1)+\dots+q_n(x_n)}\bigg)
}
is satisfied for all $(x_1,\dots,x_n)\in V$, where $V$ is an open subset of $I^n$ containing its diagonal. Consider the map $F:I^n\to J^n$ defined by $F(x_1,\dots.x_n):=(f(x_1),\dots,f(x_n)$. Then $F$ is a continuous bijection between $I^n$ and $J^n$ whose inverse is also continuous. Therefore, 
$U:=F(V)$ is an open subset of $J^n$ containing its diagonal. For $(u_1,\dots,u_n)\in U$ we have that $F^{-1}(u_1,\dots,u_n)=(f^{-1}(u_1),\dots,f^{-1}(u_1))\in V$, therefore, applying \eq{Eq1} with 
$(x_1,\dots,x_n):=(f^{-1}(u_1),\dots,f^{-1}(u_1))$, we get that
\Eq{*}{
  &\frac{(p_1\circ f^{-1})(u_1)u_1+\dots+(p_n\circ f^{-1})(u_n)u_n}{(p_1\circ f^{-1})(u_1)+\dots+(p_n\circ f^{-1})(u_n)}\\
  &=(g\circ f^{-1})^{-1}\bigg(\frac{(q_1\circ f^{-1})(u_1)(g\circ f^{-1})(u_1)+\dots+(q_n\circ f^{-1})(u_n)(g\circ f^{-1})(u_n)}{(q_1\circ f^{-1})(u_1)+\dots+(q_n\circ f^{-1})(u_n)}\bigg),
}
holds for all $(u_1,\dots,u_n)\in U$.
With the notations introduced in \eq{hPQ}, the above equality can be rewritten, for all $(u_1,\dots,u_n)\in U$, as
\Eq{Eq2}{
	\frac{P_1(u_1)u_1+\dots+P_n(u_n)u_n}{P_1(u_1)+\dots+P_n(u_n)}
	=h^{-1}\bigg(\frac{Q_1(u_1)h(u_1)+\dots+Q_n(u_n)h(u_n)}{Q_1(u_1)+\dots+Q_n(u_n)}\bigg),
}
which shows that the means $A_{\id,P}$ and $A_{h,Q}$ are locally equal on $J^n$. If, in addition, the means $A_{f,p}$ and $A_{g,q}$ are globally equal on $I^n$, i.e., the equality \eq{Eq1} holds on $V=I^n$, then $U=J^n$, which shows that the means $A_{\id,P}$ and $A_{h,Q}$ are globally equal on $J^n$.

The reversed implication can checked similarly, therefore its proof is omitted.
\end{proof}

In what follows, we deduce the first-order necessary condition for the equality of $n$-variable generalized Bajraktarević means.

\Thm{1NC}{Let $h:J\to\R$ a differentiable function with a nonvanishing first derivative, let $n\geq 2$, let $P,Q:J\to\R_+^n$, let $i\in\{1,\dots,n\}$ and assume that $P_i$ is continuous on $J$. Assume that the means $A_{\id,P}$ and $A_{h,Q}$ are locally equal on $J^n$. Then
\Eq{PQ}{
  \frac{P_i}{P_0}=\frac{Q_i}{Q_0},
}
where
\Eq{PQ0}{
  P_0:=P_1+\dots+P_n \qquad\mbox{and}\qquad
  Q_0:=Q_1+\dots+Q_n.
}}

\begin{proof} According to \lem{reg}, the continuity of $P_i$ implies the continuity of $Q_i$. On the other hand, by \thm{DB}, it follows that first-order partial derivatives
$\partial_i A_{\id,P}$ and $\partial_i A_{h,Q}$
exist at the diagonal points of $J^n$, and due to the local validity of the equality $A_{\id,P}=A_{h,Q}$, they are equal to each other, i.e., \eq{PQ} holds.
\end{proof}

Without assuming the differentiability of $h$, for fixed $i\in\{1,\dots,n\}$, $u\in J$, we have that
\Eq{*}{
  \frac{P_i(u)}{P_0(u)}=
  (\partial_iA_{h,Q})^\Delta(u)=\lim_{v\to u}
  \frac{1}{v-u}\bigg(h^{-1}\bigg(\frac{(h\cdot Q_0)(u)-(h\cdot Q_i)(u)+(h\cdot Q_i)(v)}{Q_0(u)-Q_i(u)+Q_i(v)}\bigg)-u\bigg).
}
Could one get the differentiability of $h$ from here at a.e. point with a nonvanishing derivative?

\Cor{1NC}{Let $h:J\to\R$ a differentiable function with a nonvanishing first derivative, let $n\geq 2$, let $P,Q:J\to\R_+^n$ and assume that $P$ is continuous on $J$. Assume that the means $A_{\id,P}$ and $A_{h,Q}$ are locally equal on $J^n$. Then $Q$ is also continuous on $J$ and there exists a continuous function $r:J\to\R_+$ such that
\Eq{PQr}{
  Q=rP.
}
Furthermore, the equality 
\Eq{Eq3}{
  h\bigg(\frac{P_1(u_1)u_1+\dots+P_n(u_n)u_n}{P_1(u_1)+\dots+P_n(u_n)}\bigg)
  =\frac{(rhP_1)(u_1)+\dots+(rhP_n)(u_n)}{(rP_1)(u_1)+\dots+(rP_n)(u_n)}.
}
holds locally on $J^n$ with respect to $u=(u_1,\dots,u_n)\in J^n$.}

\begin{proof} Define the functions $P_0$ and $Q_0$ by \eq{PQ0} and let $r:=Q_0/P_0$. Then, according to \thm{1NC}, for all $i\in\{1,\dots,n\}$, the equality \eq{PQ} holds and $Q_i$ is continuous.  Therefore the equation \eq{PQr} is valid and $r$ is also continuous. 

Assume that the means $A_{\id,P}$ and $A_{h,Q}$ are equal over an open subset $U\subseteq J^n$ which contains the diagonal of $J^n$. Then \eq{Eq2} is valid for all $u=(u_1,\dots,u_n)\in U$. Replacing $Q_i(u_i)$ by $r(u_i)P_i(u_i)$ and applying the function $h$ to this equality side by side, we can conclude that \eq{Eq3} is valid for all $u=(u_1,\dots,u_n)\in U$, i.e., it holds locally on $J^n$. 
\end{proof}

First we establish a sufficient condition for the equality \eq{Eq3}.

\Thm{SC}{If there exist four real constants $a,b,c,d$ such that, for all $u\in J$,
\Eq{rh}{
  r(u)=cu+d>0,\qquad h(u)=\frac{au+b}{cu+d}
}
and $P:J\to\R_+^n$ is arbitrary, then \eq{Eq3} holds for all $u=(u_1,\dots,u_n)\in J^n$.}

\begin{proof} Let $u=(u_1,\dots,u_n)\in J^n$. Then
\Eq{*}{
  &h\bigg(\frac{P_1(u_1)u_1+\dots+P_n(u_n)u_n}{P_1(u_1)+\dots+P_n(u_n)}\bigg)
  =\frac{a\dfrac{P_1(u_1)u_1+\dots+P_n(u_n)u_n}{P_1(u_1)+\dots+P_n(u_n)}+b}{c\dfrac{P_1(u_1)u_1+\dots+P_n(u_n)u_n}{P_1(u_1)+\dots+P_n(u_n)}+d}\\
  &\qquad=\frac{a(P_1(u_1)u_1+\dots+P_n(u_n)u_n)+b(P_1(u_1)+\dots+P_n(u_n))}{c(P_1(u_1)u_1+\dots+P_n(u_n)u_n)+d(P_1(u_1)+\dots+P_n(u_n))}\\
  &\qquad=\frac{(au_1+b)P_1(u_1)+\dots+(au_n+b)P_n(u_n)}{(cu_1+d)P_1(u_1)+\dots+(cu_n+d)P_n(u_n)}
  =\frac{(rhP_1)(u_1)+\dots+(rhP_n)(u_n)}
	{(rP_1)(u_1)+\dots+(rP_n)(u_n)}.
}
\end{proof} 

The first main goal of this paper is to show that if the equality \eq{Eq3} holds locally in $J^n$, then, under some regularity assumptions, the functions $r$ and $h$ are of the form \eq{rh} for some real constants $a,b,c,d$. To accomplish this goal, we are going to show that $r$ and $h$ are twice differentiable functions and 
\Eq{*}{
  r''=0\qquad\mbox{and}\qquad (rh)''=0.
}
The proof of this property will be split in several propositions.

\Prp{1st}{Let $h:J\to\R$ be continuously differentiable with a nonvanishing first derivative and let $r:J\to\R_+$. Let $P:J\to\R_+^n$, $i,j\in\{1,\dots,n\}$, $i\neq j$ and suppose that $P_i$ and $P_j$ are continuously differentiable. Assume that \eq{Eq3} holds locally with respect to $(u_1,\dots,u_n)\in J^n$. Then $r$ is continuously differentiable and 
\Eq{ij}{
  &\Big(P_i(u)\cdot(P_i(u)+P_j(v))+P_i'(u)\cdot P_j(v)\cdot(u-v)\Big)\\
  &\qquad\times \Big((h'rP_j)(v)\cdot((rP_i)(u)+(rP_j)(v))+(rP_j)'(v)\cdot(rP_i)(u)\cdot(h(v)-h(u))\Big)\\
  &=\Big(P_j(v)\cdot(P_i(u)+P_j(v))+P_j'(v)\cdot P_i(u)\cdot(v-u) \Big)\\
  &\qquad\times\Big((h'rP_i)(u)\cdot((rP_i)(u)+(rP_j)(v))+(rP_i)'(u)\cdot(rP_j)(v)\cdot(h(u)-h(v))\Big)
}
is valid locally for $(u,v)\in J^2$.}

\begin{proof} The local validity of the equality \eq{Eq3} implies that the $n$-variable means $A_{\id,P}$ and $A_{h,rP}$ are locally equal on $J^n$. According to \cor{RV}, it follows that two-variable means $A_{\id,(P_i,P_j)}$ and $A_{h,(rP_i,rP_j)}$ are locally equal on $J^2$.
In view of assertion (ii) of \lem{reg}, the continuous differentiability of $P_i$ implies the continuous differentiability of $rP_i$. Therefore, the function $r$ must be continuously differentiable. By local equality of these two-variable means, there exists an open subset $V\subseteq J^2$ containing the diagonal of $J^2$ such that, for all $(u,v)\in V$, we have
\Eq{*}{
  h\bigg(\frac{P_i(u)u+P_j(v)v}{P_i(u)+P_j(v)}\bigg)
  =\frac{(rhP_i)(u)+(rhP_j)(v)}{(rP_i)(u)+(rP_j)(v)}.
}
Computing the partial derivatives of this equation side by side with respect to $u$ and $v$, for all $(u,v)\in V$, we get
\Eq{*}{
  &h'\bigg(\frac{P_i(u)u+P_j(v)v}{P_i(u)+P_j(v)}\bigg)
  \times \frac{P_i(u)\cdot(P_i(u)+P_j(v))+P_j(v)\cdot P_i'(u)\cdot(u-v)}{(P_i(u)+P_j(v))^2}\\
  &\qquad=\frac{(h'rP_i)(u)\cdot((rP_i)(u)+(rP_j)(v))+(rP_i)'(u)\cdot(rP_j)(v)\cdot(h(u)-h(v))}{((rP_i)(u)+(rP_j)(v))^2}
}
and 
\Eq{*}{
  &h'\bigg(\frac{P_i(u)u+P_j(v)v}{P_i(u)+P_j(v)}\bigg)
  \times\frac{P_j(v)\cdot(P_i(u)+P_j(v))+P_i(u)\cdot P_j'(v)\cdot(v-u)}{(P_i(u)+P_j(v))^2}\\
  &\qquad=\frac{(h'rP_j)(v)\cdot((rP_i)(u)+(rP_j)(v))+(rP_j)'(v)\cdot(rP_i)(u)\cdot(h(v)-h(u))}{((rP_i)(u)+(rP_j)(v))^2}.
}
Multiplying the left and right hand sides of the first equality by the right and left hand sides of the second one, after cancellations, the stated equality \eq{ij} follows.
\end{proof}

\Prp{2nd}{Let $h:J\to\R$ be continuously differentiable with a nonvanishing first derivative and let $r,P_i,P_j:J\to\R_+$ be continuously differentiable (where $i,j\in\{1,\dots,n\}$ and $i\neq j$). Assume that \eq{ij} holds locally with respect to $(u,v)\in J^2$. Then $h'r^2$ is a constant function on $J$, hence $h$ is twice continuously differentiable on $J$, and the equality
\Eq{ij+}{
    \frac{r(v)-r(u)}{r(u)r(v)(u-v)}
  =&\frac{(rP_j')(v)\cdot(P_i^2)(u)+(rP_i')(u)\cdot(P_j^2)(v)}
  {(rP_i)(u)\cdot (rP_j)(v)\cdot(P_i(u)+P_j(v))}\\
  &-\frac{(rP_i)'(u)\cdot(rP_j^2)(v)+(rP_j)'(v)\cdot(rP_i^2)(u)}
  {P_i(u)\cdot P_j(v)\cdot((rP_i)(u)+(rP_j)(v))}\cdot\frac{1}{u-v}\int_v^u\frac{1}{r^2}\\
  &+\frac{(rP_j')(v)\cdot(r'P_i)(u) - (rP_i')(u)\cdot(r'P_j)(v)}{(P_i(u)+P_j(v))\cdot((rP_i)(u)+(rP_j)(v))}\cdot\int_v^u\frac{1}{r^2}
}
is valid locally for $(u,v)\in J^2$.}

\begin{proof} Let $V\subseteq J^2$ be an open set containing the diagonal of $J^2$ where \eq{ij} holds. Then, rearranging this equality, for all $(u,v)\in V$, we obtain
\Eq{*}{
  0=&P_i(u)\cdot P_j(v)\cdot(P_i(u)+P_j(v))\cdot((rP_i)(u)+(rP_j)(v))\cdot((h'r)(u)-(h'r)(v))\\
  &+(P_j'(v)\cdot(h'rP_i^2)(u)+P_i'(u)\cdot(h'rP_j^2)(v))\cdot((rP_i)(u)+(rP_j)(v))\cdot(v-u)\\
  &+((rP_i)'(u)\cdot(rP_j^2)(v)+(rP_j)'(v)\cdot(rP_i^2)(u))\cdot(P_i(u)+P_j(v))\cdot(h(u)-h(v))\\
  &+P_i(u)\cdot P_j(v)\cdot((rP_j')(v)\cdot(r'P_i)(u) - (rP_i')(u)\cdot(r'P_j)(v))\cdot(v-u)\cdot(h(u)-h(v)).
}
Dividing both sides of this equation by $P_i(u)\cdot P_j(v)\cdot(P_i(u)+P_j(v))\cdot((rP_i)(u)+(rP_j)(v))\cdot(u-v)$, we get
\Eq{hr}{
  \frac{(h'r)(u)-(h'r)(v)}{u-v}
  =&\frac{P_j'(v)\cdot(h'rP_i^2)(u)+P_i'(u)\cdot(h'rP_j^2)(v)}
  {P_i(u)\cdot P_j(v)\cdot(P_i(u)+P_j(v))}\\
  &-\frac{(rP_i)'(u)\cdot(rP_j^2)(v)+(rP_j)'(v)\cdot(rP_i^2)(u)}
  {P_i(u)\cdot P_j(v)\cdot((rP_i)(u)+(rP_j)(v))}\cdot\frac{h(u)-h(v)}{u-v}\\
  &+\frac{(rP_j')(v)\cdot(r'P_i)(u) - (rP_i')(u)\cdot(r'P_j)(v)}
  {(P_i(u)+P_j(v))\cdot((rP_i)(u)+(rP_j)(v))}\cdot(h(u)-h(v)).
}
Let $u\in J$ be fixed. Observe that the limit of the right hand side as $v\to u$ exists. (Here we use that $V$ is open and contains the point $(u,u)$.) This shows that $h'r$ is differentiable and hence $h'$ is differentiable, too. This proves that $h$ is twice differentiable. Upon taking the limit as $v\to u$, we get that the equality
\Eq{*}{
  h''r+h'r'=(h'r)'
  =&h'r\cdot\frac{P_j'P_i^2+P_i'P_j^2}
  {P_iP_j(P_i+P_j)}-h'\cdot\frac{(rP_i)'P_j^2+(rP_j)'P_i^2}{P_iP_j(P_i+P_j)}=-h'r'
}
is valid over $J$.
Therefore,
\Eq{*}{
  h''r+2h'r'=0,
}
which implies that $h'r^2$ is a constant function on $J$. Using that $r$ is continuously differentiable, it follows that $h'$ is also continuously differentiable and hence $h$ is twice continuously differentiable.

Denote $\gamma:=h'r^2\neq0$. Then, we get that $h'=\gamma/r^2$, thus the equality \eq{hr} can be rewritten as \eq{ij+} for all $(u,v)\in V$.
\end{proof}

In what follows, we are going to prove that, for $w\in I$, the \emph{symmetric derivative of $r'$ at $w$}, i.e., the limit
\Eq{*}{
   \lim_{t\to0}\frac{r'(w+t)-r'(w-t)}{2t}
}
exists and equals $0$. Then, using the quasi-mean value theorem of Aull \cite[Theorem 1]{Aul67} (see also \cite[Theorem 6.3]{SahRie98}) from the theory of symmetric derivatives, it will follow that $r'$ is a constant function.

\Prp{3rd}{Let $r,P_i,P_j:J\to\R_+$ be continuously differentiable (where $i,j\in\{1,\dots,n\}$ and $i\neq j$) such that $P_i(u)\neq P_j(u)$ for all $u\in J$. Assume that \eq{ij+} holds locally with respect to $(u,v)\in J^2$. Then $r$ is an affine function on $J$, i.e., there exists real constants $c,d$ such that $r(u)=cu+d$ holds for all $u\in J$.}

\begin{proof} Let $V\subseteq J^2$ be an open set containing the diagonal of $J^2$ where \eq{ij+} holds. Multiplying this equality by $((rP_i)(u)+(rP_j)(v))$ side by side and then rearranging it, for all $(u,v)\in V$, we arrive at the following equation
\Eq{Pr}{
   \frac{P_i(u)}{r(v)}&\bigg(\frac{r(u)-r(v)}{u-v}-r'(v)\bigg)+\frac{P_j(v)}{r(u)}\bigg(\frac{r(u)-r(v)}{u-v}-r'(u)\bigg)\\
   &=\bigg(\frac{P_i(u)r'(v)}{r(v)}+\frac{P_j(v)r'(u)}{r(u)}+\frac{P_i'(u)P_j^3(v)+P_j'(v)P_i^3(u)}{(P_i(u)+P_j(v))P_i(u)P_j(v)}\bigg)\bigg(\frac{r(u)r(v)}{u-v}\int_v^u\frac{1}{r^2}-1\bigg)\\
  &\qquad+\frac{P_i'(u)P_j(v)r(u)}{(P_i(u)+P_j(v))r(v)}\bigg(\frac{(r'(v)(u-v)+r(v))r(v)}{u-v}\int_v^u\frac{1}{r^2}-1\bigg)\\
  &\qquad+\frac{P_j'(v)P_i(u)r(v)}{(P_i(u)+P_j(v))r(u)}\bigg(\frac{(r'(u)(v-u)+r(u))r(u)}{u-v}\int_v^u\frac{1}{r^2}-1\bigg).
}
Let $w\in J$ be fixed and substitute $u:=w+t$ and $v:=w-t$ into \eq{Pr}. First, compute the limit of both sides of this equality as $t\to0$. By the continuity of the function $r$, we have that 
\Eq{uv}{
  \lim_{t\to 0}\frac{1}{2t}\int_{w-t}^{w+t}\frac{1}{r^2}=\frac{1}{r(w)^2},
}
hence both sides of the equality \eq{Pr} tend to zero as $t\to 0$. In what follows, we show that right hand side of the equation \eq{Pr} multiplied by $\frac{1}{2t}$ tends to zero as $t\to 0$. For this purpose, we compute the following three limits using \eq{uv}, the equality
\Eq{*}{
  \frac{d}{dt}\bigg(\int_{w-t}^{w+t}\frac{1}{r^2}\bigg)
  =\frac{d}{dt}\bigg(\int_{w}^{w+t}\frac{1}{r^2}-\int_{w}^{w-t}\frac{1}{r^2}\bigg)
  =\frac{1}{r(w+t)^2}+\frac{1}{r(w-t)^2},
}
and L'Hospital's Rule:
\Eq{*}{
  \lim_{t\to 0}&\frac{1}{2t}\bigg(\frac{r(w+t)r(w-t)}{2t}\int_{w-t}^{w+t}\frac{1}{r^2}-1\bigg)
  =\lim_{t\to 0}\frac{1}{4t^2}\bigg(r(w+t)r(w-t)\int_{w-t}^{w+t}\frac{1}{r^2}-2t\bigg)\\
 &=\lim_{t\to 0}\frac{1}{8t}\bigg(\big(r'(w+t)r(w-t)-r(w+t)r'(w-t)\big)\int_{w-t}^{w+t}\frac{1}{r^2}\\
 &\qquad\qquad\qquad+r(w+t)r(w-t)\bigg(\frac{1}{r(w+t)^2}+\frac{1}{r(w-t)^2}\bigg)-2\bigg)\\
 &=\frac{r'(w)}{4r(w)}-\frac{r'(w)}{4r(w)}
 +\lim_{t\to 0}\frac{1}{8t}\bigg(\frac{r(w-t)}{r(w+t)}+\frac{r(w+t)}{r(w-t)}-2\bigg)\\
 &=\lim_{t\to 0}\frac{(r(w+t)-r(w-t))^2}{8tr(w+t)r(w-t)}=\frac{1}{r^2(w)}\lim_{t\to 0}\frac{(r(w+t)-r(w-t))^2}{8t}=0.
}
In a similar fashion, we also get
\Eq{*}{
  &\lim_{t\to 0}\frac{1}{2t}\bigg(\frac{(-2tr'(w+t)+r(w+t))r(w+t)}{2t}\int_{w-t}^{w+t}\frac{1}{r^2}-1\bigg)\\
  &=-\lim_{t\to 0}\bigg(\frac{r'(w+t)r(w+t)}{2t}\int_{w-t}^{w+t}\frac{1}{r^2}\bigg)+\lim_{t\to 0}\frac{1}{4t^2}\bigg(r(w+t)^2\int_{w-t}^{w+t}\frac{1}{r^2}-2t\bigg)\\
  &=-\frac{r'(w)}{r(w)}
  +\lim_{t\to 0}\frac{1}{8t}\bigg(2r'(w+t)r(w+t)\int_{w-t}^{w+t}\frac{1}{r^2}+r(w+t)^2\bigg(\frac{1}{r(w+t)^2}+\frac{1}{r(w-t)^2}\bigg)-2\bigg)\\
  &=-\frac{r'(w)}{r(w)}+\frac{r'(w)}{2r(w)}
  +\lim_{t\to 0}\frac{r(w+t)-r(w-t)}{2t}\cdot\frac{r(w+t)+r(w-t)}{4r(w-t)^2}
  =-\frac{r'(w)}{2r(w)}+\frac{r'(w)}{2r(w)}=0
}
and
\Eq{*}{
  &\lim_{t\to 0}\frac{1}{2t} \bigg(\frac{(2tr'(w-t)+r(w-t))r(w-t)}{2t}\int_{w-t}^{w+t}\frac{1}{r^2}-1\bigg)\\
  &=\lim_{t\to 0}\bigg(\frac{r'(w-t)r(w-t)}{2t}\int_{w-t}^{w+t}\frac{1}{r^2}\bigg)+\lim_{t\to 0}\frac{1}{4t^2}\bigg(r(w-t)^2\int_{w-t}^{w+t}\frac{1}{r^2}-2t\bigg)\\
  &=\frac{r'(w)}{r(w)}
  +\lim_{t\to 0}\frac{1}{8t}\bigg(-2r'(w-t)r(w-t)\int_{w-t}^{w+t}\frac{1}{r^2}+r(w-t)^2\bigg(\frac{1}{r(w+t)^2}+\frac{1}{r(w-t)^2}\bigg)-2\bigg)\\
  &=\frac{r'(w)}{r(w)}-\frac{r'(w)}{2r(w)}
  +\lim_{t\to 0}\frac{1}{8t}\frac{r(w-t)^2-r(w+t)^2}{r(w+t)^2}
  =\frac{r'(w)}{r(w)}-\frac{r'(w)}{2r(w)}-\frac{r'(w)}{2r(w)}=0.
}

Using the above three equalities and the continuity of $P_i,P_i'$ and $r$, it follows that the right hand side of the equation \eq{Pr} at $u=w+t$ and $v:=w-t$ multiplied by $\frac{1}{2t}$ tends to zero as $t\to 0$. Therefore, the same property is valid for the left hand side of \eq{Pr}, i.e, for all $w\in J$,
\Eq{*}{
  0&=\lim_{t\to 0}\bigg(\frac{P_i(w+t)}{r(w-t)}\bigg(\frac{r(w+t)-r(w-t)}{4t^2}-\frac{r'(w-t)}{2t}\bigg)\\
  &\qquad\qquad +\frac{P_j(w-t)}{r(w+t)}\bigg(\frac{r(w+t)-r(w-t)}{4t^2}-\frac{r'(w+t)}{2t}\bigg)\bigg).
}
Replacing $t$ by $-t$ and multiplying the equality so obtained by $-1$, we also obtain that
\Eq{*}{
  0&=\lim_{t\to 0}\bigg(\frac{P_i(w-t)}{r(w+t)}\bigg(\frac{r(w+t)-r(w-t)}{4t^2}-\frac{r'(w+t)}{2t}\bigg)\\
  &\qquad\qquad+\frac{P_j(w+t)}{r(w-t)}\bigg(\frac{r(w+t)-r(w-t)}{4t^2}-\frac{r'(w-t)}{2t}\bigg)\bigg).
}
These equalities, with obvious manipulations of the limits, imply that
\Eq{*}{
  0&=\lim_{t\to 0}\bigg(\frac{P_j(w-t)}{P_i(w+t)}\frac{r(w-t)}{r(w+t)}\bigg(\frac{r(w+t)-r(w-t)}{4t^2}-\frac{r'(w+t)}{2t}\bigg)\\
  &\hspace{5.5cm}+\frac{r(w+t)-r(w-t)}{4t^2}-\frac{r'(w-t)}{2t}\bigg)
}
and
\Eq{*}{
  0&=\lim_{t\to 0}\bigg(\frac{P_i(w-t)}{P_j(w+t)}\frac{r(w-t)}{r(w+t)}\bigg(\frac{r(w+t)-r(w-t)}{4t^2}-\frac{r'(w+t)}{2t}\bigg)\\
  &\hspace{5.5cm} +\frac{r(w+t)-r(w-t)}{4t^2}-\frac{r'(w-t)}{2t}\bigg).
}
Subtracting the two equalities side by side, we arrive at
\Eq{*}{
  0=\lim_{t\to 0}\frac{r(w-t)}{r(w+t)}\bigg(\frac{P_j(w-t)}{P_i(w+t)}-\frac{P_i(w-t)}{P_j(w+t)}\bigg)\bigg(\frac{r(w+t)-r(w-t)}{4t^2}-\frac{r'(w+t)}{2t}\bigg),
}
which yields that
\Eq{*}{
  0=\lim_{t\to 0}\big(P_j(w+t)P_j(w-t)-P_i(w+t)P_i(w-t)\big)\bigg(\frac{r(w+t)-r(w-t)}{4t^2}-\frac{r'(w+t)}{2t}\bigg).
}
Replacing $t$ by $-t$ and then multiplying the equality so obtained by $-1$, we also obtain that
\Eq{*}{
  0=\lim_{t\to 0}\big(P_j(w+t)P_j(w-t)-P_i(w+t)P_i(w-t)\big)\bigg(\frac{r(w+t)-r(w-t)}{4t^2}-\frac{r'(w-t)}{2t}\bigg).
}
Subtracting the second equality from the first one, we conclude that
\Eq{w}{
   0=\lim_{t\to 0}\big(P_j(w+t)P_j(w-t)-P_i(w+t)P_i(w-t)\big)\frac{r'(w+t)-r'(w-t)}{2t}.
}

By our assumption, $P_i(w)\neq P_j(w)$ holds for all $w\in J$. Thus we have that
\Eq{*}{
  \lim_{t\to 0}\big(P_j(w+t)P_j(w-t)-P_i(w+t)P_i(w-t)\big)=P_j^2(w)-P_i^2(w)\neq0,
}
therefore, the limit of the second factor in \eq{w} exists and equals zero. This shows that the symmetric derivative of $r'$ vanishes on $J$. According the quasi-mean value theorem of Aull \cite[Theorem 1]{Aul67} (see also \cite[Theorem 6.3]{SahRie98}), this implies that $r'$ is constant on $J$, hence $r$ is an affine function on $J$.
\end{proof}

\Thm{Main}{Let $f,g:I\to\R$ be strictly monotone continuous functions. Let $n\ge2$ and $p=(p_1,\dots,p_n),q=(q_1,\dots,q_n):I\to\R_+^n$. Assume that $p$ is continuous on $I$ and there exist indices $i,j\in\{1,\dots,n\}$, $i\neq j$ and a nonempty open subinterval $I_0$ of $I$ such that
\begin{enumerate}[(a)]
 \item $p_i(x)\neq p_j(x)$ for all $x\in I_0$. 
 \item $p_i\circ f^{-1}$ and $p_j\circ f^{-1}$ are continuously differentiable on the interval $f(I_0)$.
 \item $g\circ f^{-1}$ is continuously differentiable on the interval $f(I_0)$ with a nonvanishing first derivative.
\end{enumerate}
Then the following assertions are equivalent to each other:
\begin{enumerate}[(i)]
 \item The means $A_{f,p}$ and $A_{g,q}$ are globally equal on $I^n$.
 \item The means $A_{f,p}$ and $A_{g,q}$ are locally equal on $I^n$.
 \item There exist $a,b,c,d\in\R$ with $ad\neq bc$ such that the conditions in \eq{fg} hold on $I$.
\end{enumerate}}

\begin{proof} The implication (i)$\Rightarrow$(ii) is obvious. The implication (iii)$\Rightarrow$(i) is a consequence of \thm{BM3}. Therefore, we may restrict our attention to the proof of implication (ii)$\Rightarrow$(iii).

Assume now that assertion (ii) holds, i.e., the $n$-variable means $A_{f,p}$ and $A_{g,q}$ are locally equal on $I^n$ and hence, on $I_0^n$ as well. 
Define the interval $J_0:=f(I_0)$, and the functions $h$, $P=(P_1,\dots,P_n)$, and $Q=(Q_1,\dots,Q_n)$ on the interval $J_0$ by \eq{hPQ}. Then, applying \thm{Red}, we get that the $n$-variable means $A_{\id,P}$ and $A_{h,Q}$ are locally equal on $J_0^n$. According to our assumptions (c), $h$ is continuously differentiable on $J_0$ with a nonvanishing first derivative. The continuity of $p$, implies that $P$ is continuous on $J_0$. Thus, by the first assertion of \lem{reg}, it follows that $Q$ is also continuous on $J_0$ and there exists a continuous function $r:J_0\to\R_+$ such that \eq{PQr} holds on $J_0$. By assumption (b), we have that $P_i$ and $P_j$ are continuously differentiable on $J_0$. This, applying the second assertion of \lem{reg}, implies that $Q_i$ and $Q_j$ are also continuously differentiable on $J_0$. Therefore, $r$ is also continuously differentiable on $J_0$. In addition, by assumption (a), we also have that  $P_i(u)\neq P_j(u)$ for all $u\in J_0$. 

As a consequence of \cor{RV}, we obtain that the two-variable means $A_{\id,(P_i,P_j)}$ and $A_{h,(Q_i,Q_j)}$ are locally equal on $J_0^2$. Applying \cor{1NC}, it follows that
\Eq{*}{
  h\bigg(\frac{P_i(u)u+P_j(v)}{P_i(u)+P_j(v)}\bigg)
  =\frac{(rhP_i)(u)+(rhP_j)(v)}{(rP_i)(u)+(rP_j)(v)}
}
holds locally for $(u,v)\in J_0^2$.

Using now \prp{1st}, we obtain that \eq{ij} is valid 
locally for $(u,v)\in J_0^2$. Then, in view of \prp{2nd}, we get that $h'r^2$ is a constant function, $h$ is twice differentiable on $J_0$, and \eq{ij+} is satisfied locally for $(u,v)\in J_0^2$.
Finally, using \prp{3rd}, we conclude that, for some $c,d\in\R$, the function $r$ is of the form $r(u)=cu+d$.
Using that $h'r^2$ is a constant on $J_0$, it follows that 
\Eq{*}{
  (hr)''=h''r+2h'r'+hr''=h''r+2h'r'=\frac1r(h'r^2)'=0,
}
which shows that $hr$ is also an affine function, i.e.,
there exist constants $a,b\in\R$ such that $(hr)(u)=au+b$ for all $u\in J_0$. This proves that \eq{rh} is valid on $J_0$. Therefore, with the substitution $u:=f^{-1}(x)$, using also that \eq{PQr} is valid on $J_0$, we can obtain that \eq{fg} is valid on the interval $I_0$.
Then, applying \thm{BM4}, we can see that the equality \eq{fg} is also valid on $I$. 

This completes the proof of the implication (ii)$\Rightarrow$(iii).
\end{proof}

\Thm{Main+}{Let $f,g:I\to\R$ be strictly monotone continuous functions. Let $n\ge3$ and $p=(p_1,\dots,p_n)\to\R_+^n$, $q=(q_1,\dots,q_n)\to\R_+^n$. Assume that $p$ is continuous on $I$ and there exist indices $i,j,k\in\{1,\dots,n\}$ with $i<j<k$ and a nonempty open subinterval $I_0$ of $I$ such that
\begin{enumerate}[(a)]
 \item $p_i\circ f^{-1}$, $p_j\circ f^{-1}$ and $p_k\circ f^{-1}$ are continuously differentiable on the interval $f(I_0)$.
 \item $g\circ f^{-1}$ is continuously differentiable on the interval $f(I_0)$ with a nonvanishing first derivative.
\end{enumerate}
Then the following three assertions are equivalent to each other:
\begin{enumerate}[(i)]
 \item The means $A_{f,p}$ and $A_{g,q}$ are globally equal on $I^n$.
 \item The means $A_{f,p}$ and $A_{g,q}$ are locally equal on $I^n$.
 \item There exist $a,b,c,d\in\R$ with $ad\neq bc$ such that the conditions in \eq{fg} hold on $I$.
\end{enumerate}}

\begin{proof} The implications (i)$\Rightarrow$(ii) and (iii)$\Rightarrow$(i) follow in the same way as in the proof of the previous theorem. 

Assume now that assertion (ii) holds, i.e., the $n$-variable means $A_{f,p}$ and $A_{g,q}$ are locally equal on $I^n$ and hence, on $I_0^n$ as well. 
Define the interval $J_0:=f(I_0)$, and the functions $h$, $P=(P_1,\dots,P_n)$, and $Q=(Q_1,\dots,Q_n)$ on the interval $J_0$ by \eq{hPQ}. Then, using the same argument as in the proof of the previous theorem, 
we can establish that $h$ is continuously differentiable on $J_0$ with a nonvanishing first derivative, $P$ and $Q$ are continuous on $J_0$, $P_i,P_j,P_k,Q_i,Q_j,Q_k$ are continuously differentiable and there exists a continuously differentiable function $r:J_0\to\R_+$ such that \eq{PQr} holds on $J_0$. 

Then, in view of \thm{RV}, it follows that the 3-variable means $A_{f,(p_i,p_j,p_k)}$ and $A_{g,(q_i,q_j,q_k)}$ are locally equal on $I^3$. Let $x_0\in I_0$ be fixed. The equalities
\Eq{*}{
  p_i(x_0)=p_j(x_0)+p_k(x_0),\qquad
  p_j(x_0)=p_i(x_0)+p_k(x_0),\qquad
  p_k(x_0)=p_i(x_0)+p_j(x_0)
}
cannot hold simultaneously (because the sum of the left hand sides is strictly smaller than the sum of the right hand sides). By the symmetric role of these indices, we may assume that $p_i(x_0)\neq p_j(x_0)+p_k(x_0)$. Then, by the continuity of the functions $p_i,p_j,p_k$ at $x_0$, the inequality $p_i(x)\neq p_j(x)+p_k(x)$ holds for all $x$ belonging to a neighborhood of $x_0$. Without loss of generality, we may also assume that it holds for all $x\in I_0$.

The local equality of the 3-variable means $A_{f,(p_i,p_j,p_k)}$ and $A_{g,(q_i,q_j,q_k)}$ easily implies that the 2-variable means $A_{f,(p_i,p_j+p_k)}$ and $A_{g,(q_i,q_j+q_k)}$ are locally equal on $I^2$. Observe that we can apply the previous theorem to these 2-variable means and we have that assertion (ii) holds in this setting. Thus, the assertion (iii) in this setting shows that there exist four constants $a,b,c,d\in R$ such that
\Eq{*}{
  cf+d>0,\qquad g=\frac{af+b}{cf+d},\qquad
  q_i=(cf+d)p_i,\qquad 
  q_j+q_k=(cf+d)(p_j+p_k)
}
hold on $I_0$. The equality $q_i=(cf+d)p_i$ implies that
$Q_i(u)=(cu+d)P_i(u)$ for all $u\in J_0=f(I_0)$. Therefore, $r(u)=cu+d$ for all $u\in J_0$. This equality together with \eq{PQr} yield that, for all $u\in J_0$,
\Eq{*}{
  Q(u)=(cu+d)P(u).
}
Therefore, we can deduce that
\Eq{*}{
  q=(cf+d)p
}
is valid on $I_0$ and hence \eq{fg} holds on $I_0$. Now, applying \thm{BM4}, it follows that \eq{fg} also holds on $I$, i.e., assertion (iii) of our theorem is valid.
\end{proof}

\def\MR#1{}


\providecommand{\bysame}{\leavevmode\hbox to3em{\hrulefill}\thinspace}
\providecommand{\MR}{\relax\ifhmode\unskip\space\fi MR }
\providecommand{\MRhref}[2]{%
  \href{http://www.ams.org/mathscinet-getitem?mr=#1}{#2}
}
\providecommand{\href}[2]{#2}

\end{document}